# Normal and Triangular Determinantal Representations of Multivariate Polynomials

*Massimo Salvi*

*Department of Mathematics and Computer Science*

*90123-Palermo, PA.*

*Italy*

e-mail: *massimo.salvi1@istruzione.it*

**Abstract.**

In this paper we give a new and simple algorithm to put any multivariate polynomial into a normal determinant form in which each entry has the form $a_i x_i + b_i$, and in each column the same variable appears. We also apply the algorithm to obtain a triangular determinant representation, a reduced determinant representation, and a uniform determinant representation of any multivariable polynomial. The algorithm could be useful for obtaining representations of dimensions smaller than those available up to now to solve numerical problems.

## 1.Introduction

The problem of finding a determinantal representation of polynomials by mean of linear polynomials as entries has been studied by many scholars [3,4,12,13, 14]. Recently the problem of representing polynomials in determinantal form has received an increasing interest in order to solve many problems in various areas of mathematics [5,6,7,10]. More recent results have been obtained by other authors [1,2,9] who defined the *uniform determinantal representation* of a multivariate polynomial and gave an algorithm to compute this representation having the aim to find the minimum possible matrix dimension. An application of this representation has been found in order to solve systems of multivariate polynomials [8, 1]. The problem of finding the minimal dimension of any given multivariate polynomial is still open, in this paper we give some examples that show the potential of the algorithm in order to find determinantal representation of small dimension.

In this paper we prove this main result:

*Any multivariate polynomial with integer coefficients can be written as determinant of a matrix whose entries are binomials of the form $a_i x_i + b_i$, where we denote with $x_i$ the variables of the polynomial, and $a_i, b_i$ are integers. Moreover in each column of this matrix appear binomials in the same variable.* In this paper we define this form as *normal determinantal representation* of the polynomial.

We prove this result by presenting a simple algorithm that permits the costruction of such a matrix. Starting from the previous normal determinal representation we define another determinantal representation that we define *triangular determinant representation* of the polynomial. We generalize these representations to multivariate polynomials with coeffients in a field. Then we show how to reduce a polynomial into uniform determinantal representation by using its trianguar determinantal form and we conjecture that the dimension of the obtained matrix is the minimum.

The paper is organized as follows:

In section 2 we give a fundamental result on integer matrices which is the core of the proposed algorithm to reduce any multivariate polinomials into normal determinantal representation.

In section 3 we define the *chain-form* and the *improved chain-form* of a multivariate polynomial which are used in the algorithm.

In section 4 we define the *normal determinantal representation* of a polynomial and the algorithm to reduce a multivariate polynomials with integer coefficients into normal determinantal representation is presented.

In section 5 we define the *triangular determinantal representation* of a polynomial and show how to obtain such a representation starting from the normal determinantal representation.

In section 6 we generalize the previous results and the agorithm to multivariate polynomials with coefficients in a field.

In section 7 we show how to reduce the dimension of the matrix and obtain a reduced determinantal representation starting from the triangular one.

In section 8 we show how to obtain a uniform determinantal representation starting from the reduced determinantal representation.

In section 9 we conclude with some open questions and a conjecture about the dimension of the matrix obtained for the reduced determinantal representation.

## 2. Some preliminary results

In this section we will recall some results that the author has already shown elsewhere [11] and play an important role in the algorithm that will be shown in what follows. The results will be presented in a slightly different form from the original work and will be proven for completeness.

If A is a square matrix we denote with $|A|$ the determinant of A.

**Theorem 1**. The diophantine equation

$$\begin{vmatrix} b_{11} & b_{12} & & b_{1n} \\ b_{21} & b_{22} & & .. \\ .. & .. & & \\ b_{n-11} & b_{n-12} & & b_{n-1n} \\ a_1 & a_2 & & a_n \end{vmatrix} = 1 \qquad (1)$$

where the entries $a_i$ are integers and the $b_{ij}$ are the unknowns, is solvable in Z if and only if $\gcd\{a_1, a_2.., a_n\} = 1$.

*Proof.* If $d = \gcd\{a_1, a_2.., a_n\}$ is greater than 1 then the determinant would be dividible by d and it's in contradiction with its value equal to 1, so the condition is clearly necessary.

Let's find explicitly a solution in the case that $\gcd\{a_1, a_2.., a_n\} = 1$.

Since $\gcd\{a_1, a_2.., a_n\} = 1$, it's possible (see Appendix 1), by applying the euclidean algorithm on the set $\{a_1, a_2.., a_n\}$ to compute the $\gcd\{a_1, a_2.., a_n\}$, to find an unimodular integer matrix $M$ such that:

$(a_1 \quad a_2 \quad .. \quad .. \quad a_n)M = (0 \quad 0 \quad .. \quad .. \quad 1)$

The matrix $M$ is unimodular, it means that its determinant is +1 or -1, if the determinant of $M$ is -1 we can always change the signs of the first column in order to obtain a matrix $M$ with determinant

+1 without changing the validity of equation 1b (since we change the signs of a column whose product with the vector $(a_1 \quad a_2 \quad .. \quad .. \quad a_n)$ gives 0 ). This fact will be usefull in what follows.

If we take into cosideration the inverse $M^{-1}$, which is unimodular and integer too, we have

$$(a_1 \quad a_2 \quad .. \quad .. \quad a_n) = (0 \quad 0 \quad .. \quad .. \quad 1)M^{-1}$$

and we can build a solution of the equation (*1) by considering the product

$$IM^{-1} = \begin{bmatrix} 1 & 0 & .. & .. & 0 \\ 0 & 1 & & & .. \\ .. & & & & .. \\ .. & & & 1 & 0 \\ 0 & 0 & .. & 0 & 1 \end{bmatrix} M^{-1} = \begin{bmatrix} JM^{-1} \\ (0 \quad 0 \quad .. \quad .. \quad 1)M^{-1} \end{bmatrix} = \begin{bmatrix} JM^{-1} \\ (a_1 \quad a_2 \quad .. \quad .. \quad a_n) \end{bmatrix} \quad (1b)$$

Moreover, we have that $|M|=1$, hence $|M^{-1}|=1$, therefore $|IM^{-1}|=|I|\cdot|M^{-1}|=1$, hence the matrix $JM^{-1}$ give a solution of (*1). We notice that the matrix $JM^{-1}$ is simply the matrix of the first $n$-1 rows of $M^{-1}$ ∎

**Theorem 2.** Given the linear form $a_1x_1 + a_2x_2 + .. + a_nx_n$, in which all values $\{a_1, a_2.., a_n\}$ are integers, it is always possible to find an integer matrix $A$, $(n-1) \times n$, such that:

$$a_1x_1 + a_2x_2 + .. + a_nx_n = \begin{vmatrix} a_{11} & .. & .. & .. & a_{1n} \\ . & & & & . \\ . & & & & . \\ a_{n-11} & .. & .. & .. & a_{n-1n} \\ x_1 & x_2 & .. & .. & x_n \end{vmatrix} = \begin{vmatrix} A \\ X \end{vmatrix} \quad (2)$$

*Proof.* We can suppose that the $\gcd\{a_1,a_2..,a_n\}=1$, in fact given any form $a'_1 x_1 + a'_2 x_2 + .. + a'_n x_n$ and setting $k = \gcd\{a'_1, a'_2.., a'_n\}$ we can always consider the new form: $\dfrac{a'_1 x_1 + a'_2 x_2 + .. + a'_n x_n}{k} = a_1x_1 + a_2x_2 + .. + a_nx_n$ in which $\gcd\{a_1, a_2.., a_n\}=1$. If we suppose that a matrix $A$ exists such that the previous (2) holds, then by multiplying the first row (or any other) of the matrix $A$ by $k$, we obtain that (2) holds for $a'_1 x_1 + a'_2 x_2 + .. + a'_n x_n$ as well.

Now let's consider the matrix $M$ defined as follows:

$$M = \begin{bmatrix} 1 & 0 & .. & .. & 0 \\ 0 & 1 & 0 & .. & 0 \\ .. & & 0 & .. & .. & 0 \\ 0 & & .. & 0 & 1 & 0 \\ M_{n1} & M_{n2} & .. & .. & M_{nn} \end{bmatrix}$$

in the last row the entries are:
$$\begin{cases} M_{n1} = b_{11}x_1 + b_{12}x_2 + .. + b_{1n}x_n \\ M_{n2} = b_{21}x_1 + b_{22}x_2 + .. + b_{2n}x_n \\ \quad .. \\ M_{nn-1} = b_{n-11}x_1 + b_{n-12}x_2 + .. + b_{n-1n}x_n \\ M_{nn} = a_1x_1 + a_2x_2 + .. + a_nx_n \end{cases}$$

where $b_{ij}$ are integer values which will be found in what follows,

the rest of the entries are: $\begin{cases} M_{ij} = 1 & \text{if} \quad i = j \\ M_{ij} = 0 & \text{if} \quad i \neq j \end{cases}$

By definition we have:
$$\begin{vmatrix} 1 & 0 & .. & .. & 0 \\ 0 & 1 & 0 & .. & 0 \\ .. & & 0 & .. & .. & 0 \\ 0 & & .. & 0 & 1 & 0 \\ M_{n1} & M_{n2} & .. & .. & M_{nn} \end{vmatrix} = M_{nn} = a_1x_1 + a_2x_2 + .. + a_nx_n$$

hence the theorem will be proved if we are able to find an unimodular matrix $U$:

$$U = \begin{bmatrix} u_{11} & u_{12} & .. & .. & u_{1n} \\ u_{21} & .. & .. & .. & .. \\ .. & .. & .. & .. & .. \\ .. & .. & .. & .. & .. \\ u_{n1} & u_{n2} & .. & .. & u_{nn} \end{bmatrix} \text{ such that } |U| = 1 \text{ and: } \begin{bmatrix} 1 & 0 & .. & .. & 0 \\ 0 & 1 & 0 & .. & 0 \\ .. & & 0 & .. & .. & 0 \\ 0 & & .. & 0 & 1 & 0 \\ M_{n1} & M_{n2} & .. & .. & M_{nn} \end{bmatrix} U = \begin{bmatrix} a_{11} & a_{12} & .. & .. & a_{1n} \\ a_{21} & .. & .. & .. & .. \\ .. & .. & .. & .. & .. \\ a_{n-11} & .. & .. & .. & a_{n-1n} \\ x_1 & x_2 & .. & .. & x_n \end{bmatrix}$$

as far as the last row is concerned, we have:

$$[b_1x_1+b_{12}x_2+..+b_{1n}x_n; \quad .. \quad .. \quad b_{n-11}x_1+b_{n-12}x_2+..+b_{n-1n}x_n; \quad a_1x_1+a_2x_2+..+a_nx_n] \begin{bmatrix} u_{11} & u_{12} & .. & .. & u_{1n} \\ u_{21} & .. & .. & .. & .. \\ .. & .. & .. & .. & .. \\ .. & .. & .. & .. & .. \\ u_{n1} & u_{n2} & .. & .. & u_{nn} \end{bmatrix} = [x_1 \quad x_2 \quad .. \quad .. \quad x_n]$$

by computing this product we obtain:

$$\begin{cases} (b_1u_{11}+b_2u_{21}+..+b_{n-11}u_{n-11}+a_1u_{n1})x_1 + (b_{12}u_{11}+b_{22}u_{21}+..+b_{n-12}u_{n-11}+a_2u_{n1})x_2 +..+ (b_{1n}u_{11}+b_{2n}u_{21}+..+b_{n-1n}u_{n-11}+a_nu_{n1})x_n = x_1 \\ (b_1u_{12}+b_2u_{22}+..+b_{n-11}u_{n-12}+a_1u_{n2})x_1 + (b_{12}u_{12}+b_{22}u_{22}+..+b_{n-12}u_{n-12}+a_2u_{n2})x_2 +..+ (b_{1n}u_{12}+b_{2n}u_{22}+..+b_{n-1n}u_{n-12}+a_nu_{n2})x_n = x_2 \\ \quad\quad .. \\ \quad\quad .. \\ (b_1u_{1n}+b_2u_{2n}+..+b_{n-11}u_{n-1n}+a_1u_{nn})x_1 + (b_{12}u_{1n}+b_{22}u_{2n}+..+b_{n-12}u_{n-1n}+a_2u_{nn})x_2 +..+ (b_{1n}u_{1n}+b_{2n}u_{2n}+..+b_{n-1n}u_{n-1n}+a_nu_{nn})x_n = x_n \end{cases}$$

we can verify that the previous system is equivalent to requiring that:

$$\begin{bmatrix} u_{11} & u_{21} & .. & .. & u_{n1} \\ u_{12} & .. & .. & & u_{n2} \\ .. & & & & .. \\ .. & .. & .. & .. & .. \\ u_{1n} & & & & u_{nn} \end{bmatrix} \begin{bmatrix} b_{11} & b_{12} & & b_{1n} \\ b_{21} & b_{22} & & .. \\ .. & .. & & \\ b_{n-11} & b_{n-12} & & b_{n-1n} \\ a_1 & a_2 & & a_n \end{bmatrix} = \begin{bmatrix} 1 & 0 & .. & 0 & 0 \\ 0 & 1 & & & \\ .. & & .. & & \\ .. & & & .. & 0 \\ 0 & & & 0 & 1 \end{bmatrix} = I \qquad (3)$$

or 
$$U^T \begin{bmatrix} b_{11} & b_{12} & & b_{1n} \\ b_{21} & b_{22} & & .. \\ .. & .. & & \\ b_{n-11} & b_{n-12} & & b_{n-1n} \\ a_1 & a_2 & & a_n \end{bmatrix} = I$$

Now, since we required that $|U|=1$, it follows that $|U^T|=1$, and by using Binet's theorem in (3) we have:

$$\begin{vmatrix} b_{11} & b_{12} & & b_{1n} \\ b_{21} & b_{22} & & .. \\ .. & .. & & \\ b_{n-11} & b_{n-12} & & b_{n-1n} \\ a_1 & a_2 & & a_n \end{vmatrix} = 1 \qquad (4)$$

But in the equation (4) $\gcd\{a_1, a_2.., a_n\}=1$, so we can apply Theorem 1 and conclude that (4) is solvable. Namely, there exists the required set of integer entries $b_{ij}$ and, therefore, there exists the required unimodular matrix $U$ as well. In fact from (3) we have:

$$(U^T)^{-1} = \begin{bmatrix} b_{11} & b_{12} & & b_{1n} \\ b_{21} & b_{22} & & .. \\ .. & .. & & \\ b_{n-11} & b_{n-12} & & b_{n-1n} \\ a_1 & a_2 & & a_n \end{bmatrix} \qquad (5)$$

and the inverse of an unimodular matrix is unimodular as well, hence $U^T$ is unimodular, therefore $U$ is unimodular as well. ∎

Now we are ready to prove the following theorem that plays a fundamental role in our algorithm for representing a polynomial in determinantal form.

**Theorem 3**. In order to obtain the matrix $A$, defined as in Theorem 2, we first compute the matrix $M$, as described in the proof of Theorem 1, then the matrix $A$ is given by the first $n$-1 rows of $M^T$.

*Proof.* We can find explicitly the matrix $A$ :

first, considering the previous equation
$$\begin{bmatrix} 1 & 0 & .. & .. & 0 \\ 0 & 1 & 0 & .. & 0 \\ .. & 0 & .. & .. & 0 \\ 0 & .. & 0 & 1 & 0 \\ M_{n1} & M_{n2} & .. & .. & M_{nn} \end{bmatrix} U = \begin{bmatrix} a_{11} & a_{12} & .. & .. & a_{1n} \\ a_{21} & .. & .. & .. & .. \\ .. & .. & .. & .. & .. \\ a_{n-11} & .. & .. & .. & a_{n-1n} \\ x_1 & x_2 & .. & .. & x_n \end{bmatrix}$$

, we see that $A$ is given by the first $n$-1 rows of the matrix $U$. From equation (5) we can write

$$U = \left( \begin{bmatrix} b_{11} & b_{12} & & b_{1n} \\ b_{21} & b_{22} & & .. \\ .. & .. & & \\ b_{n-11} & b_{n-12} & & b_{n-1n} \\ a_1 & a_2 & & a_n \end{bmatrix}^{-1} \right)^T$$

from (1b) we have that
$$\begin{bmatrix} b_{11} & b_{12} & & b_{1n} \\ b_{21} & b_{22} & & .. \\ .. & .. & & \\ b_{n-11} & b_{n-12} & & b_{n-1n} \\ a_1 & a_2 & & a_n \end{bmatrix} = M^{-1}$$

, so we have simply: $U = ((M^{-1})^{-1})^T = M^T$, therefore $A$ is obtained by the first n-1 rows of U ∎

Let's do an example to clarify this result:

**Example 2.1.** If we consider $2x_1 - 7x_2 + 4x_3$, by using the generalized euclidean algorithm (Appendix 1) we can find the matrix $M$ with determinant equal to 1, it is $\begin{bmatrix} -12 & 7 & -3 \\ -4 & 2 & -1 \\ -1 & 0 & 0 \end{bmatrix}$, in fact we have:

$$[2 \; -7 \; 4] \begin{bmatrix} -12 & 7 & -3 \\ -4 & 2 & -1 \\ -1 & 0 & 0 \end{bmatrix} = [0 \; 0 \; 1]$$

Hence considering the first 2 rows of $M^T$ we obtain the requested matrix $A$:

$$A = \begin{bmatrix} -12 & -4 & -1 \\ 7 & 2 & 0 \end{bmatrix}$$

and we can verify that: $\begin{vmatrix} -12 & -4 & -1 \\ 7 & 2 & 0 \\ x_1 & x_2 & x_3 \end{vmatrix} = 2x_1 - 7x_2 + 4x_3$

**Corollary 1.** Be $P(x_1 \; .. \; x_k) = \sum_{i=1}^{l} c_i x_0^{n_{0i}} x_1^{n_{1i}} ... x_k^{n_{ki}}$ a general polynomial in $k$ variables, with $l$ monomials and $c_i \in Z$, then it is always possible to find an integer matrix $A$, $(l-1) \times l$, such that:

$$P(x_1 \; .. \; x_k) = \begin{vmatrix} a_{11} & .. & .. & .. & a_{1l} \\ . & & & & . \\ . & & & & . \\ a_{l-11} & .. & .. & .. & a_{l-1l} \\ x_0^{n_{01}} x_1^{n_{11}} ... x_k^{n_{k1}} & x_0^{n_{02}} x_1^{n_{12}} ... x_k^{n_{k2}} & .. & .. & x_0^{n_{0l}} x_1^{n_{1l}} ... x_k^{n_{kl}} \end{vmatrix}$$

## 3. Chains of monomials and chain-form of polynomials

**Definition 1.** Given the monomial $m = c \, x_1^{n1} x_2^{n2} .. x_k^{nk}$ we define *chain of the monomial m* an ordered set of monomials starting from $x_1^{n1} x_2^{n2} .. x_k^{nk}$ and ending to $x_k$ in which each monomial is obtained by the previous (starting from the left) by dividing for one of the variables.

If a monomial is an integer costant c then the chain is simply c.

**Example 3.1.** If we are given the monomial $2x_1^2 x_2^3 x_3^2$, a chain could be the ordered set

$$(x_1^2 x_2^3 x_3^2 \quad x_1^1 x_2^3 x_3^2 \quad x_2^3 x_3^2 \quad x_2^2 x_3^2 \quad x_2^1 x_3^2 \quad x_3^2 \quad x_3^1)$$

It's clear that for every monomial is possible to find a chain.

Notice also that given a monomial m in n variables and one of the variables in m, say $x_i$, it's always possible to find a chain of m ending with $x_i$ (the proof is trivial).

**Definition 2.** If we are given a general polynomial in $k$ variables and $l$ terms $P(x_1 \ .. \ x_k) = \sum_{i=1}^{l} c_i x_0^{n_{0i}} x_1^{n_{1i}} ... x_k^{n_{ki}}$ with integer coefficients $c_i$, we define the **chain-form of** $P(x_1 \ .. \ x_k)$ as follows:

1. We group the terms with the same degree and consider these groups by descendent order of degree.
2. Starting from the monomials of maximum degree we find, for each monomial $m_j$, a chain of $m_j$.
3. The chain-form of $P(x_1 \ .. \ x_k)$ is formed by summing the monomials of these chains taken with coefficient 0 if the corresponding term does not appear in $P(x_1 \ .. \ x_k)$ and with coefficient $c_i$ if it appears in $P(x_1 \ .. \ x_k)$ with coefficient and such a term will not be taken into consideration in the next steps.
4. We consider the monomials of lower degree and go on like in step 2, without considering the mononials with coefficient different from 0 eventually found in step 3.

**Example 3.1.** Be $P(x_1 \ x_2) = 3x_1^3 x_2^2 - 4x_1^2 x_2^3 + x_1^2 x_2^2 - 5x_1 x_2^2 + 2x_1^3 + 2x_1 x_2 + 2$

| Degree | Groups | Chains | Terms of the chain-form of $P(x_1 \ x_2)$ |
|---|---|---|---|
| 5 | $3x_1^3 x_2^2$ | C1: $x_1^3 x_2^2 \quad x_1^2 x_2^2 \quad x_1^1 x_2^2 \quad x_2^2 \quad x_2^1$ | $3x_1^3 x_2^2 \ +1x_1^2 x_2^2 \ +0x_1^1 x_2^2 \ +0x_2^2 \ +0x_2^1$ |
|   | $-4x_1^2 x_2^3$ | C2: $x_1^2 x_2^3 \quad x_1^1 x_2^3 \quad x_2^3 \quad x_2^2 \quad x_2^1$ | $-4x_1^2 x_2^3 \ +0x_1^1 x_2^3 \ +0x_2^3 \ +0x_2^2 \ +0x_2^1$ |
| 4 | $x_1^2 x_2^2$ | present in C1 | |
| 3 | $-5x_1 x_2^2$ | present in C1 | |
|   | $+2x_1^3$ | C3: $x_1^3 \quad x_1^2 \quad x_1^1$ | $+2x_1^3 \ +0x_1^2 \ +0x_1^1$ |
| 2 | $+2x_1 x_2$ | C4: $x_1^1 x_2^1 \quad x_2^1$ | $+2x_1^1 x_2^1 \ +0x_2^1$ |
| 0 | $+2$ | C5: $+2$ | $+2$ |

Hence we can write a chain-form of the considered polynomial as:

$P(x_1 \ x_2) = 3x_1^3 x_2^2 \ +1x_1^2 x_2^2 \ +0x_1^1 x_2^2 \ +0x_2^2 \ +0x_2^1 -4x_1^2 x_2^3 \ +0x_1^1 x_2^3 \ +0x_2^3 \ +0x_2^2 \ +0x_2^1 +$
$+2x_1^3 \ +0x_1^2 \ +0x_1^1 \ +2x_1^1 x_2^1 \ +0x_2^1 \ +2$

**Definition 3.** If we are given a general polynomial in $k$ variables and $l$ terms $P(x_1 \; .. \; x_k) = \sum_{i=1}^{l} c_i x_0^{n_{0i}} x_1^{n_{1i}} ... x_k^{n_{ki}}$ with integer coefficients $c_i$, then we define the *improved chain-form* of $P(x_1 \; .. \; x_k)$ as follows:

1. We group the terms with the same degree and consider these groups by descendent order of degree.
2. Starting from the monomials of maximum degree we find, for each monomial $m_j$, a chain $C_j$ of $m_j$. If a monomial $m_k$ in the chain of $m_j$ already appears in the chain $C_i$ of a previous monomial $m_i$, then we stop the chain of $m_j$ to $m_k$ and build a unique chain by using $C_i$ and $C_j$ the two as follows: we consider the monomials of $C_i$ from the first to the one preceding $m_k$, then we add the momomials of $C_j$ from the first to $m_k$ and, finally, we add the monomials of $C_i$ from the one that follows $m_k$ to the last.
3. The chain-form of $P(x_1 \; .. \; x_k)$ is formed by summing the monomials of these chains taken with coefficient 0 if the corresponding term does not appear in $P(x_1 \; .. \; x_k)$ and with coefficient $c_i$ if it appears in $P(x_1 \; .. \; x_k)$ with coefficient and such a term will not be taken into consideration in the next step.
4. We consider the monomials of lower degree and go on like in step 2, without considering the mononials with coefficient different from 0 eventually found in step 3.

**Example 3.2.** Be $P(x_1 \; x_2) = 3x_1^3 x_2^2 - 4x_1^2 x_2^3 + x_1^2 x_2^2 - 5x_1 x_2^2 + 2x_1^3 + 2x_1 x_2 + 2$

| Degree | Groups | Chains | Terms of the chain-form of $P(x_1 \; x_2)$ |
|---|---|---|---|
| 5 | $3x_1^3 x_2^2$ $-4x_1^2 x_2^3$ | C1: $x_1^3 x_2^2$ $x_1^2 x_2^2$ $x_1^1 x_2^2$ $x_2^2$ $x_2^1$<br>C2: $x_1^2 x_2^3$ $x_1^2 x_2^2$ $x_1^1 x_2^2$ $x_2^2$ $x_2^1$<br>In step2 we find that $x_1^2 x_2^2$ appears in C1, so we build a unique chain :<br>C1+2:<br>$x_1^3 x_2^2$ $x_1^2 x_2^3$ $x_1^2 x_2^2$ $x_1^1 x_2^2$ $x_2^2$ $x_2^1$ | $3x_1^3 x_2^2 - 4x_1^2 x_2^3 + 1x_1^2 x_2^2 - 5x_1^1 x_2^2 + 0x_2^2 + 0x_2^1$ |
| 4 | $x_1^2 x_2^2$ | appears in C1 | |
| 3 | $-5x_1 x_2^2$ $+2x_1^3$ | appears in C1<br>C3: $x_1^3$ $x_1^2$ $x_1^1$ | $+2x_1^3 + 0x_1^2 + 0x_1^1$ |
| 2 | $2x_1 x_2$ | C4: $x_1^1 x_2^1$ $x_2^1$<br>In step 2 we find that $x_2^1$ appears in C1+2<br>so we build a unique chain: | |

| | | C1+2+4 $x_1^3x_2^2$ $x_1^2x_2^3$ $x_1^2x_2^2$ $x_1^1x_2^2$ $x_1^1x_2^1$ $x_2^2$ $x_2^1$ | $3x_1^3x_2^2 - 4x_1^2x_2^3 + x_1^2x_2^2 - 5x_1^1x_2^2 + x_1^1x_2^1$ $+ 0x_2^2 + 0x_2^1$ |
|---|---|---|---|
| 0 | +2 | C5: +2 | +2 |

Hence we can write an improved chain-form of the considered polynomial as:

$$P(x_1; x_2) = 3x_1^3x_2^2 - 4x_1^2x_2^3 + x_1^2x_2^2 - 5x_1^1x_2^2 + x_1^1x_2^1 + 0x_2^2 + 0x_2^1 + 2x_1^3 + 0x_1^2 + 0x_1^1 + 2$$

From step 2 it follows that, in general, an improved chain-form is shorter than the chain-form it derives from. In example 3.1 the chain-form contains 16 monomials, in example 3.2 the improved chain-form contains 11 monomials. Finding an improved chain-form of a given polynomial is usefull in order to reduce the dimension of the matrix used to find its determinantal representation. Notice that it is possible to find different chain-forms (or improved chain-forms) for a given polynomial, but this does not affect what follows.

## 4. Normal Determinantal Representation of a multivariate polynomial

We are ready to present an algorithm which permits to build a determinantal form of any multivariate polynomial. Given a polynomial $P(x_1 \ .. \ x_k)$, we will build a square matrix with entries of the form $a_{ij}x_k + b_{ij}$ whose determinant is equal to $P(x_1 \ .. \ x_k)$. In order to do this let's consider a chain-form of $P(x_1 \ .. \ x_k)$ and let's denote by C the vector of the ordered set of the coefficients and by M the vector of the ordered set of the monomials regardless to the coefficients they are multiplied by of the chain-form as they appears at the end of the algorithm explained in §4 applied to the given $P(x_1 \ .. \ x_k)$.

For example, if we consider the following chain-form polynomial

$$P(x_1 \ x_2) = 3x_1^3x_2^2 + 1x_1^2x_2^2 + 0x_1^1x_2^2 + 0x_2^2 + 0x_2^1 - 4x_1^2x_2^3 + 0x_1^1x_2^3 + 0x_2^3 + 0x_2^2 + 0x_2^1 +$$
$$+ 2x_1^3 + 0x_1^2 + 0x_1^1 + 2x_1^1x_2^1 + 0x_2^1 + 2$$

the vectors C and M are :

$C = (3;\ 1;\ 0;\ 0;\ 0;-4;\ 0;\ 0;\ 0;\ 0;\ 2;\ 0;\ 0;\ 2;\ 0\ +2)$

$M = (x_1^3x_2^2;\ x_1^2x_2^2;\ x_1^1x_2^2;\ x_2^2;\ x_2^1;\ x_1^2x_2^3;\ x_1^1x_2^3;\ x_2^3;\ x_2^2;\ x_2^1;\ x_1^3\ x_1^2;\ x_1^1;\ x_1^1x_2^1;\ x_2^1\ 1)$

Notice that the costant c in the polynomial gives rise to a c costant in the vector C and a monomial equal to 1 in the vector M; the reason will be clear in what follows.

Using the result proven in §2 we can find a matrix A of dimension (n-1)n, where n is the dimension of the vector C, such that :

$$P(x_0,..,x_k) = \begin{vmatrix} A \\ M \end{vmatrix}$$

Let's write the esplicitly the entries and denote by $m_i$ the monomials of the vector M:

$$P(x_0,..,x_k) = \begin{vmatrix} a_{11} & .. & .. & .. & a_{1n} \\ . & & & & . \\ . & & & & . \\ a_{n-11} & .. & .. & .. & a_{n-1n} \\ m_1 & m_2 & .. & .. & m_n \end{vmatrix}$$

The ordered set of monomials $m_i$ (the vector M) is a sequence of chains, it means that, starting from the first column of the matrix, if $m_1$ has a degree greater than 1, it exists a variable $x_i$ such that $m_1 = x_i m_2$, hence, by operating the column transformation $c_1 = c_1 - x_i c_2$, (which doesn't change the value of the determinant) the matrix becomes:

$$\begin{vmatrix} a_{11} - x_i a_{12} & .. & .. & .. & a_{1n} \\ a_{21} - x_i a_{22} & & & & . \\ . & & & & . \\ a_{n-11} - x_i a_{n-11} & .. & .. & .. & a_{n-1n} \\ 0 & m_2 & .. & .. & m_n \end{vmatrix}$$

in this matrix the first column contains only binomials of the form $a_i x_i + b_i$ with the same variable $x_i$, and the determinant has not changed.

Now we can proceed in the same way with the successive columns by moving to the right and obtain the requested form of the matrix. Note that if the monomial $m_1$ has degree 1 or 0 (i.e is a costant) the column doesn't change, but it contains only integer numbers and $m_1$ itself. The final matrix will have the form:

$$\begin{bmatrix} a_{11} - x_i a_{12} & a_{12} - x_j a_{13} & a_{1k} & .. & a_{1n} \\ a_{21} - x_i a_{22} & a_{22} - x_j a_{23} & a_{2k} & & . \\ . & & & & \\ a_{n-11} - x_i a_{n-12} & a_{n-12} - x_j a_{n-13} & .. & .. & a_{n-1n} \\ 0 & 0 & x_z & & x_k \end{bmatrix}$$

Therefore we have proved, by direct costruction, the following:

**Theorem 3.** Any multivariate polynomial with integer coefficients can be written as determinant of a matrix whose entries are binomials of the form $a_i x_i + b_i$, where we denote with $x_i$ the variables of the polynomial, and $a_i, b_i$ are integers. Moreover in each column of this matrix appear binomials with the same variable (which can be different for each column).

We define this representation *normal determinantal representation* (NDR) of the polynomial.

**Example4.1**. If $P(x_1 \quad x_2) = x_1^2 + 2x_1x_2 + x_2^2$ a chain-form- of $P$ is :

$$P(x_1 \quad x_2) = x_1^2 + 0x_1 + 2x_1^1 x_2^2 + 0x_2 + 1x_2^2 + 0x_2$$

and by applying the described algorithm we obtain the $6 \times 6$ matrix:

$$\begin{bmatrix} -1 & 0 & 0 & 0 & 1 & 0 \\ -x_1 & 1 & 0 & 0 & 0 & 0 \\ 0 & 0 & 1 & 0 & -2 & 0 \\ 0 & 0 & -x_1 & 1 & 0 & 0 \\ 0 & 0 & 0 & 0 & -x_2 & 1 \\ 0 & x_1 & 0 & x_2 & 0 & x_2 \end{bmatrix}$$

In order to reduce the dimension of the matrix we can apply the previous algorithm to the improved chain-form of the polynomial, in this case only a change occurs; in the column transformation $c_n = c_n - x_i c_m$ the columns are not cecessarily adjacent. It follows from the costruction of an improved chain-form (Definition 3) that, for any monomial $m_i$ of degree greater than 1 that appears in the improved chain-form, we can find a monomial $m_j$ at the right side of $m_i$ in the chain and a variable $x_k$ such that $m_i = x_k m_j$. It means that the column transformation used in the algorithm $c_n = c_n - x_i c_m$ will use two not necessarily adjacent columns.

**Example4.2.** We can find an improved chain-form for te polynomial of Example5.1, it's :

$$P(x_1 \quad x_2) = x_1^2 + 2x_1 x_2 + 0x_1 + x_2^2 + 0x_2$$

and by applying the described algorithm we obtain the $5 \times 5$ matrix:

$$\begin{bmatrix} 1 & 0 & 0 & -1 & 0 \\ 0 & 0 & 1 & -2 & 0 \\ -x_1 & 1 & -x_2 & 0 & 0 \\ 0 & 0 & 0 & -x_2 & 1 \\ 0 & x_1 & 0 & 0 & x_2 \end{bmatrix}$$

**Example4.3.** If we consider the following polynomial:

$$P(x_1 \quad x_2) = 3x_1^3 x_2^2 - 4x_1^2 x_2^3 + x_1^2 x_2^2 - 5x_1 x_2^2 + 2x_1^3 + 2x_1 x_2$$

a possible improved chain-form is :

$$P(x_1 \quad x_2) = 3x_1^3 x_2^2 - 4x_1^2 x_2^3 + x_1^2 x_2^2 - 5x_1 x_2^2 + 2x_2 x_2 + 2x_1^3 + 0x_1^2 + 0x_1$$

by applying the described algorithm we obtain the $8 \times 8$ matrix, in which we have reordered and grouped the columns by variable:

$$\begin{bmatrix} -3x_1-1 & 3 & 0 & 0 & 0 & -3x_2 & 0 & 0 \\ 4x_1 & -4 & 0 & 0 & 0 & 4x_2-1 & 0 & 0 \\ 0 & 0 & 0 & -x_1 & 1 & 0 & 0 & -x_2 \\ 5x_1 & x_1-5 & 0 & 0 & 0 & 5x_2 & -1 & 0 \\ 2x_1 & -2 & 0 & 0 & 0 & 2x_2 & -x_2 & 1 \\ 2x_1 & -2 & 1 & 0 & 0 & 2x_2 & 0 & 0 \\ 0 & 0 & -x_1 & 1 & 0 & 0 & 0 & 0 \\ 0 & 0 & 0 & 0 & x_1 & 0 & 0 & 0 \end{bmatrix}$$

## 5. *Triangular Determinantal Representation* of a multivariate polynomial

In this section we will give an algorithm to put the normal determinantal representation of a polynomial with integer coefficients into another representation which we define *triangular determinantal representation*.

**Definition 4**: given a normal determinantal representation $N$ of a polynomials with integer coefficients, we say that it's a *triangular determinantal representation* (TDR) if:

1. it exists an integer $k \geq 0$ such that $b_{ii} \neq 0 \ \forall i, \ i \leq k$, where $b_{ii}$ are the coefficients of the variable in the entries on the diagonal. Moreover the rows $j$ such that $j > k$, contain only numbers;

2. in the column $i$, with $i \leq k$, (starting from left) the entries $N(j;i)$, where $j > i$, dont' contain the variable, i.e. are numbers;

3. in the column $i$, with $1 < i \leq k$, (starting from left) the coefficients $b_{ji}$ of the variable in the entries $N(j;i)$, where $j < i$, satisfies : $-|b_{ii}| < b_{ji} < |b_{ii}|$.

This means that the representation is in the form :

$$\begin{bmatrix} a_{11}+x_lb_{11} & a_{12}+x_jb_{12} & a_{13}+x_lb_{13} & a_{1n-1}+x_wb_{1n-1} & a_{1n}+x_zb_{1n} \\ a_{21} & a_{22}+x_jb_{22} & a_{23}+x_lb_{23} & & . \\ . & a_{32} & a_{33}+x_lb_{33} & & . \\ a_{n-11} & . & .. & a_{n-1n-1} & a_{n-1n-1} \\ a_{n1} & a_{n2} & .. & a_{nn-1} & a_{nn} \end{bmatrix}$$

We prove the following:

**Theorem 4**. Any multivariate polynomial with integer coefficients can be written in a *triangular determinantal representation* .

We prove Theorem 4 by showing an algorithm that permits to build a triangular determinantal representation of any multivariate polynomial.

*Proof.* Be $P(x_1 \;..\; x_k)$ the polynomial and $N$ a NDR (normal determinantal representation) of $P$.

We can set $N = \begin{bmatrix} a_{11} + x_i b_{12} & a_{12} + x_j b_{13} & a_{1k} & .. & a_{1n} \\ a_{21} + x_i b_{22} & a_{22} + x_j b_{23} & a_{2k} & & . \\ . & & & & \\ a_{n-11} + x_i b_{n-12} & a_{n-12} + x_j b_{n-13} & .. & .. & a_{n-1n} \\ 0 & 0 & & x_z & x_k \end{bmatrix}$

and introduce the matrices A and B:

$A = \begin{bmatrix} a_{11} & a_{12} & a_{1k} & .. & a_{1n} \\ a_{21} & a_{22} & a_{2k} & & . \\ . & & & & \\ a_{n-11} & a_{n-12} & .. & .. & a_{n-1n} \\ 0 & 0 & 0 & & 0 \end{bmatrix} \quad B = \begin{bmatrix} x_i b_{11} & x_j b_{12} & 0 & .. & 0 \\ x_i b_{11} & x_j b_{22} & 0 & & . \\ . & & & & \\ x_i b_{n-11} & x_j b_{n-12} & .. & .. & 0 \\ 0 & 0 & & x_z & x_k \end{bmatrix}$

so that we can write : $N = A + B$

$\begin{bmatrix} a_{11} + x_i b_{11} & a_{12} + x_j b_{12} & a_{1k} & .. & a_{1n} \\ a_{21} + x_i b_{21} & a_{22} + x_j b_{22} & a_{2k} & & . \\ . & & & & \\ a_{n-11} + x_i b_{n-11} & a_{n-12} + x_j b_{n-12} & .. & .. & a_{n-1n} \\ 0 & 0 & & x_z & x_k \end{bmatrix} = \begin{bmatrix} a_{11} & a_{12} & a_{1k} & .. & a_{1n} \\ a_{21} & a_{22} & a_{2k} & & . \\ . & & & & \\ a_{n-11} & a_{n-12} & .. & .. & a_{n-1n} \\ 0 & 0 & 0 & & 0 \end{bmatrix} + \begin{bmatrix} x_i b_{11} & x_j b_{12} & 0 & .. & 0 \\ x_i b_{21} & x_j b_{22} & 0 & & . \\ . & & & & \\ x_i b_{n-11} & x_j b_{n-12} & .. & .. & 0 \\ 0 & 0 & & x_z & x_k \end{bmatrix}$

We now describe an algorithm to put the polynomial $P$ in TDR (triangular determinantal representation) starting from the matrix M.

*Algorithm:*

**Step 1-** Let's consider the matrix *B*, and start with *i*=1, starting from the row number *i*, we seek the entry of the row *i* which is different from 0 and is the minimum in absolute value, let's denote with m the position in the row of such entry; two situations are possibile:

> 1.1. if the entries of the considered row of *B* are all 0, it means that the row of N contains only costants and we move it to the bottom of N and shift the others of one position toward the top, then we go to step 1, in order to do these operations on N we can operate in the same way on both matrices *A* and *B*;
>
> 1.2. if the entries of the considered row of *B* are not all 0, then the miminum exists and we swap the columns that contains the minimum and the column number *i*, in other words we swap the columns *i* and m (for these operations too, we can swap the same columns in both matrices *A* and *B*);

we apply the generalized euclidean algorithm (Appendix 1) to the entries $b_{ji}$ of the column number *i* of the matrix *B* in which $j \geq i$. By using the generalized euclidean algorithm we find the gcd which will be in the position *n* in the column *i*. If necessary we swap the rows number i and number n in order to put the gcd on the diagonal of the matrix B. After having used the generalized euclidean algorithm, the entries $b_{ji}$ of the column number *i*, such that $j > i$, will be equal to 0. As it's known, all previous operations can be obtained by multiplying the matrix *B* by an unimodular matrix *W* such that:

$$WB = W \begin{bmatrix} x_i b_{11} & x_j b_{12} & 0 & .. & 0 \\ x_i b_{21} & x_j b_{22} & 0 & & . \\ . & & & & \\ x_i b_{n-11} & x_j b_{n-12} & .. & .. & 0 \\ 0 & 0 & x_z & & x_k \end{bmatrix} = \begin{bmatrix} x_i \gcd(b_{11}, b_{21}..., b_{n-1,1}) & x_j b''_{13} & 0 & .. & 0 \\ 0 & x_i \gcd(b_{22}..., b_{n-1,2}) & 0 & & . \\ . & 0 & & b_{ll} x_z & \\ 0 & . & & 0 & .. & 0 \\ 0 & 0 & & 0 & & b_{ss} x_k \end{bmatrix}$$

The operations have to be abbplied also on the matrix A, by multiplying this matrix by W. At the end of the algorithm the request 2 of the definition 4 is satisfied.

**Step 2-** For each entry $b_{ik}$ in the column number *i* that satisfies $k < i$, we consider the rest of the division by $b_{ii}$ which is different from 0 :

$b_{ik} = q_{ik}b_{ii} + r_{ik}$ with $q_{ik}, r_{ik} \in Z$ and $-|b_{ii}| < r_{ik} < |b_{ii}|$

If we apply to the rows $R_k$ of the matrix B the transformations $R_k = R_k - q_{ik}R_i$, we obtain a matrix of the requested form by the condition 3 of the definition 4. Also in this case the operations can be obtained by multiplying the matrix $B$ by an unimodular matrix $W'$ such that:

$$W'B = W'\begin{bmatrix} x_i \gcd(b_{11}, b_{21}..., b_{n-1,1}) & x_j b''_{12} & 0 & .. & 0 \\ 0 & x_j \gcd(b_{22}..., b_{n-1,2}) & 0 & & . \\ . & 0 & b_{ll}x_z & & = \\ 0 & . & 0 & .. & 0 \\ 0 & 0 & 0 & & b_{ss}x_w \end{bmatrix} =$$

$$= \begin{bmatrix} x_i \gcd(b_{11}, b_{21}..., b_{n-1,1}) & x_j b'''_{12} & x_k b'''_{13} & .. & 0 \\ 0 & x_j \gcd(b_{22}..., b_{n-1,2}) & x_k b'''_{23} & & . \\ . & 0 & x_k \gcd(b_{33}..., b_{n-1,3}) & & \\ 0 & . & 0 & .. & 0 \\ 0 & 0 & 0 & & b_{ss}x_w \end{bmatrix}$$

that satisfies:

$-|\gcd(b_{22}..., b_{n-1,2})| < b'''_{12} < |\gcd(b_{22}..., b_{n-1,2})|$

$-|\gcd(b_{33}..., b_{n-1,3})| < b'''_{13} < |\gcd(b_{33}..., b_{n-1,3})|$, $-|\gcd(b_{33}..., b_{n-1,3})| < b'''_{23} < |\gcd(b_{33}..., b_{n-1,3})|$

and so on.

Naturally, we have to apply on A the same operations we applied on $B$.

**Step 3-** If all the entries of the rows of B number k with $k > i$ are 0, then the algorithm ends, otherwise we set $i=i+1$ and go to step 1.

When the algoritm ends we can eventually multiply the first column of the matrices $A$ and $B$ in order to obtain the right sign of the determinant (since it could have been changed in consequence of the column an d row permutations ). Finally we can find an unimodular matrix $M$ such that if we consider the product $MN = M(A + B) = MA + MB$ we obtain:

$$MN = \begin{bmatrix} a'_{11} & a'_{12} & .. & .. & a'_{1n} \\ a'_{21} & a'_{22} & .. & & . \\ . & & & & \\ a'_{n-11} & a'_{n-12} & .. & .. & a'_{n-1n} \\ a'_{n1} & a'_{n2} & .. & & a'_{nn} \end{bmatrix} + \begin{bmatrix} x_i \gcd(b_{11}, b_{21}..., b_{n-1,1}) & x_j b''_{12} & 0 & .. & 0 \\ 0 & x_i \gcd(b_{22}..., b_{n-1,2}) & 0 & & . \\ . & & 0 & b_{ll} x_z & \\ 0 & & . & 0 & .. & 0 \\ 0 & 0 & 0 & & b_{nn} x_k \end{bmatrix}$$

$$= \begin{bmatrix} a'_{11}+b'_{11} x_i & a'_{12}+b'_{12} x_i & .. & & a'_{1n}+b'_{1n} x_k \\ a'_{21} & a'_{22}+b'_{22} x_i & .. & & . \\ . & & a'_{ll}+b'_{ll} x_z & & . \\ . & . & . & .. & . \\ a'_{n1} & a'_{n2} & .. & a'_{nn-1} & a'_{nn}+b'_{nn} x_k \end{bmatrix}$$

in which in the entries $N(j;i)$, where $j<i$, satisfies: $-|b'_{ii}| < b'_{ji} < |b'_{ii}|$, therefore it is a triangular determinantal representation (TDR) of $P$ ∎

**Example 5.1**. Let's consider again the polynomial of Example 4.3:

$$P(x_1, x_2) = 3x_1^3 x_2^2 - 4x_1^2 x_2^3 + x_1^2 x_2^2 - 5x_1 x_2^2 + 2x_1^3 + 2x_1 x_2$$

By applying the previous algorithm we find the TDR of $P$ as :

$$\begin{bmatrix} -x_1-1 & 1 & 1 & 0 & 0 & 0 & -x_2 & 0 \\ -5 & x_1 & 5 & 0 & 0 & -1 & 0 & 0 \\ 0 & 0 & x_1 & -1 & 0 & 0 & 0 & 0 \\ 0 & 0 & 0 & x_1 & -1 & 0 & 0 & x_2 \\ 0 & 0 & 0 & 0 & x_1 & 0 & 0 & 1 \\ 2 & 0 & -2 & 0 & 0 & x_2 & 0 & -1 \\ -2 & 0 & 3 & 0 & 0 & 0 & 0 & 0 \\ -4 & 0 & 4 & 0 & 0 & 0 & -1 & 0 \end{bmatrix}$$

# 6. Normal and Triangular Determinantal Representation of multivariate polynomials with coefficients in a field

If we want to find a normal or triangular determinantal representation of a polynomial in which the momomials are multiplied for general coefficients in a field F, say R, we can use the same algorithm explained in §5 and §6 with a little modification. It's what we explain in what follows.

Be the polynomial $P(x_1 \;..\; x_k) = \sum_{i=1}^{l} c_i x_0^{n_{0i}} x_1^{n_{1i}} ... x_k^{n_{ki}}$ with $c_i \in F$, and let's consider the polynomial

$P'(x_1 \;..\; x_k .. x_{k+1} ... x_{k+i}) = \sum_{i=1}^{l} x_0^{n_{0i}} x_1^{n_{1i}} ... x_k^{n_{ki}} x_{k+i}$ ; if we apply the previous algorithm to put the polynomial $P'$ in normal determinantal representation we will find a matrix in which the new variables $x_{k+i}$ (one for each $c_i$) will appear only in the last row (this happens if all the chains end with the new variables $x_{k+i}$) and the determinant of the matrix is $P'$, then we can set $x_{k+i} = c_i$ and obtain a normal determinantal representation of the original polynomial P. It means that in such a way we can find a normal determinantal representation in which we can use also parameters as coefficients of the polynomial and not only numerical values. We can proceed similarly to obtain the triangular determinantal representation of $P$: first we put $P'$ in TDR, then set $x_{k+i} = c_i$.

**Example 6.1** Let's consider the polinomyal $P(x_1 \; x_2) = C_{21} x_1^2 x_2 + C_{12} x_1 x_2^2 + C_{11} x_1 x_2 + C_{10} x_1$ ; by using the described algorithm we can put it in the following determinantal representation in which the $C_{ij}$ can be real numbers.

$$\begin{bmatrix} -x_1 & 0 & 0 & 0 & 1 & 1 \\ 0 & x_1 & 0 & 0 & 0 & 0 \\ 0 & -1 & x_2 & 0 & 0 & C_{10} \\ 0 & 0 & -1 & x_2 & -C_{11} & 0 \\ -1 & 0 & 0 & 0 & -C_{21} & 0 \\ 0 & 0 & 0 & -1 & -C_{12} & 0 \end{bmatrix}$$

# 7. From Triangular to Reduced Determinantal Representation

If in the TDR of a multivariate polynomial there are k (with $k \geq 1$) rows in which only monomial of zero degree appear (i.e. there are only numerical entries), we can do a further step in order to reduce the dimension of the matrix used to represent the polynomial in determinantal representation. We refer to these rows as *costant rows* of the TDR.

Be the matrix:

$$\begin{bmatrix} a'_{11}+b'_{11}x_i & a'_{12}+b'_{12}x_i & .. & .. & a'_{1n}+b'_{1n}x_k \\ a'_{21} & a'_{22}+b'_{22}x_i & .. & & . \\ . & & a'_{ll}+b'_{ll}x_z & & . \\ a'_{n-11} & a'_{n-12} & . & .. & a'_{n-1n} \\ a'_{n1} & a'_{n2} & .. & a'_{nn-1} & a'_{nn} \end{bmatrix}$$

Then we apply the euclidean algorithm already explained in Theorem 1 to the costant rows in order to put the matrix in the form:

$$\begin{bmatrix} a''_{11}+b''_{11}x_i+c_{11}x_j+.. & a''_{12}+b''_{12}x_i+c_{12}x_j+.. & .. & .. & a''_{1n}+b''_{1n}x_i+c_{1n}x_j+. \\ a''_{21}+b''_{21}x_i+c_{21}x_j+.. & a''_{22}+b''_{22}x_i+c_{22}x_j+.. & .. & & . \\ . & & a'_{ll}+b'_{ll}x_z & & . \\ 0 & 0 & . & a''_{n-1n-1} & a''_{n-1n}+b''_{n-1n}x_i+c_{n-1n}x_j+ \\ 0 & 0 & .. & 0 & a''_{nn} \end{bmatrix}$$

The last matrix has the block form $\begin{bmatrix} L & L' \\ 0 & D \end{bmatrix}$ where L and D are square matrices, the entries of L and L' are linear forms of the variables $x_i$, 0 is a matrix of all zero entries and D is an upper triangular matrix. Therefore the determinant of such a matrix is given by $\det B = \det L \cdot \det D = |L| \cdot \prod_{i=n-m}^{n} a''_{ii}$. If we denote with $n$ the dimension of the matrix that gives the TDR of the polynomial $P$ and set $d = \prod_{i=n-m}^{n} a''_{ii}$, by multiplying the first row of the matrix $L$ by $d$ we obtain a new matrix of dimension $n-k$ with determinant equal to the polynomial $P$. We define this representation of the polynomial *reduced determinantal representation* of $P$, or *RDR*.

**Example 7.1**. Let's consider again the polynomial of Example 4.3:

$$P(x_1 \quad x_2) = 3x_1^3 x_2^2 - 4x_1^2 x_2^3 + x_1^2 x_2^2 - 5x_1 x_2^2 + 2x_1^3 + 2x_1 x_2$$

As we have seen in example 5.1, the TDR of $P$ is a $8 \times 8$ matrix with two costant rows, so $k=2$. By applying the previous algorithm we can represent the polynomial with the following $6 \times 6$ matrix that gives a RDR of $P$:

$$\begin{bmatrix} 3x_1 - 4x_2 + 1 & 1 & 0 & 0 & 0 & 0 \\ 5 & x_1 & 0 & 0 & 0 & 1 \\ -2x_1 & 0 & 0 & -1 & 0 & 0 \\ 0 & 0 & -x_2 & x_1 & -1 & 0 \\ 0 & 0 & 0 & 0 & x_1 & 0 \\ -2 & 0 & 1 & 0 & 0 & -x_2 \end{bmatrix}$$

## 8. From Reduced Determinantal Representation to Uniform Determiantal Representation

We can use the previous algorithm to obtain the Uniform Determinantal Representation (UDR) of a multivariate polynomial, a *uniform determinantal representation* of a multivariate polynomial is defined in [1] as follows:

"*Consider an n-variate polynomial of degree at most d:* $p = \sum_{|\alpha| \le d} c_\alpha x^\alpha$ , *where* $x = (x_1, x_2, \ldots x_n)$, $\alpha \in Z_{\ge 0}^n$, $|\alpha| = \sum_i \alpha_i$, $x^\alpha = \prod_i x_i^{\alpha_i}$ *and where each coefficient* $c_\alpha$ *is taken from a ground field K. A determinantal representation of p is an N ×N-matrix M of the form*

$$M = A_0 + \sum_{x=1}^{n} x_i A_i \quad \text{where each } A_i \in K^{N \times N}, \text{ with det(M)} = p.\text{"}$$

We show how to obtain a UDR of a multivariate polynomial by developing an example with the general tri-variate polynomial of degree 2:

$$p(x_1, x_2, x_3) = c_{200}x_1^2 + c_{020}x_2^2 + c_{002}x_3^2 + c_{110}x_1x_2 + c_{101}x_1x_3 + c_{011}x_2x_3 + c_{100}x_1 + c_{010}x_2 + c_{001}x_3 + c_{000}$$

where we suppose that $c_{ijk} \in R$

We can define new variables and do the substitutions on *p*:

$x_4 = c_{200}x_1$; $x_5 = c_{020}x_2$; $x_6 = c_{002}x_3$; $x_7 = c_{110}x_1$; $x_8 = c_{101}x_1$; $x_9 = c_{011}x_2$; $x_{10} = c_{100}x_1$; $x_{11} = c_{010}x_2$; $x_{12} = c_{001}x_3$; $x_{13} = c_{000}$

by applying the substitutions we obtain a new polynomial $p'$ in 13 variables with degree 2:

$$p'(x_1, x_2, \ldots x_{13}) = x_4x_1 + x_5x_2 + x_6x_3 + x_7x_2 + x_8x_3 + x_9x_3 + x_{10} + x_{11} + x_{12} + x_{13}$$

A possible improved chain-form of $p'$ is:

$$p'(x_1, x_2, \ldots x_{12}) = x_4 x_1 + 0 x_1 + x_5 x_2 + x_7 x_2 + 0 x_2 + x_6 x_3 + x_8 x_3 + x_9 x_3 + 0 x_3 + x_{10} + x_{11} + x_{12} + x_{13}$$

Now we can find a RDR of $p'$ as explained in §7 and obtain the matrix:

$$\begin{bmatrix} 0 & x_4 & 0 & -1 \\ -1 & x_5 + x_7 & 0 & 0 \\ 0 & x_6 + x_8 + x_9 & 1 & 0 \\ x_2 & x_{10} + x_{11} + x_{12} + x_{13} & -x_3 & x_1 \end{bmatrix}$$

so we can write the UDR of $p'$ as :

$$\begin{bmatrix} 0 & c_{200} x_1 & 0 & -1 \\ -1 & c_{110} x_1 + c_{020} x_2 & 0 & 0 \\ 0 & c_{101} x_1 + c_{011} x_2 + c_{002} x_3 & 1 & 0 \\ x_2 & c_{100} x_1 + c_{010} x_2 + c_{001} x_3 + c_{000} & -x_3 & x_1 \end{bmatrix}$$

**Example 8.1**

Let's consider the general bi-variate polynomial of degree 4:

$$p(x_1, x_2) = c_{40} x_1^4 + c_{31} x_1^3 x_2 + c_{22} x_1^2 x_2^2 + c_{13} x_1 x_2^3 + c_{04} x_2^4 + c_{30} x_1^3 + c_{12} x_1 x_2^2 + c_{21} x_1^2 x_2 + c_{03} x_2^3 +$$

$$+ c_{20} x_1^2 + c_{11} x_1 x_2 + c_{02} x_2^2 + c_{10} x_1 + c_{01} x_2 + c_{00}$$

we define new variables and do the substitutions on $p$:

$x_3 = c_{40} x_1 \quad x_4 = c_{31} x_1 \quad x_5 = c_{22} x_1 \quad x_6 = c_{13} x_1 \quad x_7 = c_{04} x_2 \quad x_8 = c_{30} x_1 \quad x_9 = c_{21} x_1 \quad x_{10} = c_{12} x_1 \quad x_{11} = c_{03} x_2$
$x_{12} = c_{20} x_1 \quad x_{13} = c_{11} x_1 \quad x_{14} = c_{02} x_2 \quad x_{15} = c_{10} x_1 \quad x_{16} = c_{01} x_2 \quad x_{17} = c_{00}$

by applying the substitutions we obtain a new polynomial $p'$ in 17 variables with degree 4:

$$p'(x_1, x_2, \ldots x_{17}) = x_3 x_1^3 + x_4 x_1^2 x_2 + x_5 x_1 x_2^2 + x_6 x_2^3 + x_7 x_2^3 + x_8 x_1^2 + x_9 x_1 x_2 + x_{10} x_2^2 + x_{11} x_2^2 + x_{12} x_1 + x_{13} x_2 +$$

$$+ x_{14} x_2 + x_{15} + x_{16} + x_{17}$$

A possible improved chain-form is:

$$p'(x_1, x_2, \ldots x_{17}) = x_3 x_1^3 + x_4 x_1^2 x_2 + 0 x_4 x_1^2 + 0 x_8 x_1^2 + 0 x_1^2 + x_{12} x_1 + 0 x_1 + x_5 x_1 x_2^2 + 0 x_5 x_2^2 + x_6 x_2^3 + x_7 x_2^3 + 0 x_2^3 +$$

$$+ x_{10} x_2^2 + x_{11} x_2^2 + 0 x_2^2 + x_9 x_1 x_2 + 0 x_1 x_2 + x_{13} x_2 + x_{14} x_2 + 0 x_2 + x_{15} + x_{16} + x_{17}$$

we can find a RDR of $p'$ and obtain the matrix:

$$\begin{bmatrix}
-x_2 & -1 & 0 & 0 & 0 & 0 & 0 & 0 & 0 \\
x_3+x_8 & -x_4 & -1 & 0 & 0 & 0 & 0 & 0 & 0 \\
x_{12} & 0 & x_1 & 0 & 0 & 0 & 0 & 0 & -1 \\
x_1 & 0 & 0 & 0 & 0 & -1 & 0 & 0 & 0 \\
x_6+x_7 & 0 & 0 & 0 & 1 & 0 & 0 & 0 & 0 \\
x_{10}+x_{11} & 0 & 0 & 0 & -x_2 & x_5 & 0 & -1 & 0 \\
x_9 & 0 & 0 & 0 & 0 & 0 & 1 & 0 & 0 \\
x_{13}+x_{14} & 0 & 0 & 1 & 0 & 0 & -x_1 & x_2 & 0 \\
x_{15}+x_{16}+x_{17} & 0 & 0 & -x_2 & 0 & 0 & 0 & 0 & x_1
\end{bmatrix}$$

so we can write the UDR of $p'$ as :

$$\begin{bmatrix}
-x_2 & -1 & 0 & 0 & 0 & 0 & 0 & 0 & 0 \\
(c_{40}+c_{30})x_1 & -c_{31}x_1 & -1 & 0 & 0 & 0 & 0 & 0 & 0 \\
c_{20}x_1 & 0 & x_1 & 0 & 0 & 0 & 0 & 0 & -1 \\
x_1 & 0 & 0 & 0 & 0 & -1 & 0 & 0 & 0 \\
c_{13}x_1+c_{04}x_2 & 0 & 0 & 0 & 1 & 0 & 0 & 0 & 0 \\
c_{12}x_1+c_{03}x_2 & 0 & 0 & 0 & -x_2 & c_{22}x_1 & 0 & -1 & 0 \\
c_{21}x_1 & 0 & 0 & 0 & 0 & 0 & 1 & 0 & 0 \\
c_{11}x_1+c_{02}x_2 & 0 & 0 & 1 & 0 & 0 & -x_1 & x_2 & 0 \\
c_{10}x_1+c_{01}x_2+c_{00} & 0 & 0 & -x_2 & 0 & 0 & 0 & 0 & x_1
\end{bmatrix}$$

**Example 8.2**

If we consider a specific polynomial we can obtain a matrix of a dimension lower than the one obtained by considering a general polynomial.

Let's consider the polynomial in 5 variables and degree 4 (see [1], pp. 438)

$$p(x_1,x_2,x_3,x_4,x_5) = 3x_1^2 x_2 x_3 + 4x_1 x_2 x_3 + 5x_2^2 x_4 + 6x_2 x_3 x_4 + 7x_3 x_4 + 8x_5^4 + 2$$

By using the algorithm we can find the following RDR that gives also the UDR:

$$\begin{bmatrix} 0 & -1 & 0 & 0 & 0 & 0 & -3x_1-4 & 0 & 0 & 0 & 0 & 0 \\ 0 & x_3 & -1 & 0 & 0 & 0 & 0 & 0 & 0 & 0 & 0 & 0 \\ 0 & 0 & x_2 & 1 & 0 & 0 & 0 & 0 & 0 & 0 & 0 & 0 \\ 0 & 0 & 0 & 0 & 1 & 0 & -5x_4 & 0 & 0 & 0 & 0 & 0 \\ 0 & 0 & 0 & 0 & 0 & -1 & -6x_4 & 0 & 0 & 0 & 0 & 0 \\ 0 & 0 & 0 & 0 & -x_2 & x_3 & 0 & 0 & 0 & 0 & 0 & -1 \\ 0 & 0 & 0 & 0 & 0 & 0 & -7x_4 & 1 & 0 & 0 & 0 & 0 \\ 0 & 0 & 0 & 0 & 0 & 0 & -8x_5 & 0 & -1 & 0 & 0 & 0 \\ 0 & 0 & 0 & 0 & 0 & 0 & 0 & 0 & x_5 & -1 & 0 & 0 \\ 0 & 0 & 0 & 0 & 0 & 0 & 0 & 0 & 0 & x_5 & -1 & 0 \\ -1 & 0 & 0 & 0 & 0 & 0 & 0 & 0 & 0 & 0 & x_5 & 0 \\ x_5 & 0 & 0 & -x_1 & 0 & 0 & -2 & -x_3 & 0 & 0 & 0 & x_2 \end{bmatrix}$$

This example shows that by using our algorithm we can have an improvement on the dimension, since we have found a $12 \times 12$ matrix instead of a $29 \times 29$ one as obtained in [1] for the same polynomial where the authors have applied an algorithm designed for finding a determinantal representation of a general 5-variate polynomial with degree 4.

## 9. Conclusion

The algorithm we have introduced and applied in this paper is fast and simple; due also to the fact that only two matrices are needed (as explained in §5 ) to handle the matrix that gives the various different determinantal representations and only matrix operations are used. The examples proposed in the paper have been obtained with a C++ implementation of the algorithms explained. We defined some new kind of determinantal representations (normal, triangular and reduced determinantal representations) explained and obtained by a new algorithm. We think that the algorithm could be useful to reduce the dimension of the matrices used for solving numerical problems, as shown in [1] and [8], and we gave some examples in §8. Some questions remain open:

- how to find the chain-form (§3) with minimal dimension for a given polynomial?
- how to find the minimal reduced determinantal representation of a given polynomial?
- which is the relation between the lower dimension of the representation that it's possible to find with our algorithm and the lower bound of a uniform determinantal representation found in [1] ?

In our research we found that, given a polynomial, it is possible to find two reduced determinantal representations of different dimensions if we start from two improved chain-form of different dimensions, therefore it seems to exist a relation between the dimension of the chain-form of the polynomial and the dimension of its reduced determinantal representation.

We conjecture the following:

- *the algorithm explained in §7 gives rise to the minimal reduced determinantal representation if the chain-form of the polynomial is minimal*

*Appendix 1*

If the row matrix $C = \begin{pmatrix} a_1 & a_2 & ... & ... & a_n \end{pmatrix}$ is given, in which the entries are integer, in order to find the gcd of the set of integer we propose the following algorithm:

1. find the minimum, in absolute value, among the integers different from zero, be it in the position m;
2. starting from the first not zero integer and not set in the position m, compute the rest of the division of the integer by the minimum;
3. if the rest is equal to 0 substitute the integer with 0 and go on with th next entry;
4. if the rest is not equal to 0 substitute the integer with the rest and it is the new value of minimum, its position is the new value for m, go back to step 2;
5. the operations from 2 to 4 are repeated untill only one entry is different from 0;
6. the only integer different from zero is set in the last position.

It's easy to prove the following [15] :

**Theorem**. At the end of the previous algorithm, the integer in the last position is the greater common divisor (gcd) of the initals entries $a_i$ of C.

Each one of the operations in 2,3,4 and 6 can be obtained by right-multiplying C by a unimodular matrix $U_i$ which represents the trasformation. When the previous algorithm ends, by multiplying all the $U_i$ matrices we obtain a unimodular matrix $M$, $M = U_1 U_2 \cdot ... \cdot U_n$ such that:
$\begin{pmatrix} a_1 & a_2 & .. & .. & a_n \end{pmatrix} M = \begin{pmatrix} 0 & 0 & ... & \gcd(a_1 & a_2 & .. & .. & a_n) \end{pmatrix}$